\documentclass[fleqn]{article}
\usepackage{amssymb,latexsym}
\newtheorem{prop}{Proposition}[section]
\newtheorem{cor}[prop]{Corollary}
\newtheorem{theo}[prop]{Theorem}
\newtheorem{lemma}[prop]{Lemma}

\newtheorem{conj}[prop]{Conjecture}

\newtheorem{ex}[prop]{Example}
\newtheorem{note}{Remark}
\title{Sets of generators and chains of subspaces}
\author{Antonio Pasini}

\begin{document}

\maketitle

\begin{abstract}
The rank of a point-line geometry $\Gamma$ is usually defined as the generating rank of $\Gamma$, namely the minimal cardinality of a generating set. However, when the subspace lattice of $\Gamma$ satisfies the Exchange Property we can also try a different definition: consider all chains of subspaces of $\Gamma$ and take the least upper bound of their lengths as the rank of $\Gamma$. If $\Gamma$ is finitely generated then these two definitions yield the same number. On the other hand, as we shall show in this paper, if infinitely many points are needed to generate $\Gamma$ then the rank as defined in the latter way is often (perhaps always) larger than the generating rank. So, if we like to keep the first definition we should accordingly discard the second one or modify it. We can modify it as follows: consider only well ordered chains instead of arbitrary chains. As we shall prove, the least upper bound of the lengths of well ordered chains of subspaces is indeed equal to the generating rank. According to this result, the (possibly infinite) rank of a polar space can be characterized as the least upper bound of the lengths of well ordered chains of singular subspaces; referring to arbitrary chains would be an error.  
\end{abstract}

\section{Introduction}

\subsection{Definitions and results}\label{sec 1}

Following Shult \cite{Shult}, we define a {\em point-line geometry} as a pair $\Gamma = (P, {\cal L})$ where $P$ (the set of {\em points}) is a non-empty set and $\cal L$ (the set of {\em lines}) is a family of subsets of $P$, each of which contains at least two points (see also Buekenhout and Cohen \cite{BC}, where point-line geometries are called {\em line-spaces}). 

A {\em subspace} of $\Gamma$ is a subset $S\subseteq P$ such that, if a line $\ell\in{\cal L}$ meets $S$ in at least two points, then $\ell\subseteq S$. For $X\subseteq P$, we denote by $\langle X\rangle$ the subspace of $\Gamma$ {\em generated} by $X$, namely the smallest subspace of $\Gamma$ containing $X$, also calling it the {\em span} of $X$.      

Let $\mathrm{Gen}(\Gamma) = \{X\subseteq P~|~ \langle X\rangle = P\}$ be the family of sets of generators of $\Gamma$ ({\em generating sets} for short). The {\em generating rank} of $\Gamma$, hencforth denoted by $\mathrm{rk}_{\mathrm{gen}}(\Gamma)$, is the least cardinality $|X|$ of a generating set $X$: 
\[\mathrm{rk}_{\mathrm{gen}}(\Gamma) ~ := ~ \mathrm{min}(|X|~|~X\in \mathrm{Gen}(\Gamma)).\]
A subset $X\subseteq P$ is {\em independent} if $\langle Y\rangle \subset \langle X\rangle$ for every proper subset $Y\subset X$. We denote by $\mathrm{Ind}(\Gamma)$ the family of independent sets of $\Gamma$. 

The set $\mathrm{B}(\Gamma) := \mathrm{Ind}(\Gamma)\cap\mathrm{Gen}(\Gamma)$ is the family of minimal members of $\mathrm{Gen}(\Gamma)$. We call them {\em bases} of $\Gamma$. Clearly, every finite generating set contains a basis. So, if $\Gamma$ is finitely generated then $B(\Gamma) \neq \emptyset$, but geometries also exist that do not admit any basis (see Example \ref{ex1}, Subsection \ref{sec ex}).  

Let $X$ and $Y$ be two generating sets of $\Gamma$. Every $x\in X$ belongs to the span $\langle Y_x\rangle$ of a suitable finite subset $Y_x\subseteq Y$. Accordingly, the set $Y' := \cup_{x\in X}Y_x$ spans $\Gamma$ and $Y' = Y$ if $Y\in \mathrm{B}(\Gamma)$. If $X$ is finite then $Y'$ is finite as well. On the other hand, if $X$ is infinite then $|Y'| \leq \sum_{x\in X}|Y_x| = |X|$. Therefore, if $\mathrm{B}(\Gamma)\neq \emptyset$, either all bases of $\Gamma$ are finite (but possibly of different size) and every generating set contains a basis or no finite set of points generates $\Gamma$ and all bases of $\Gamma$ have the same cardinality. Hence,
\begin{equation}\label{intro 1}
\mathrm{rk}_{\mathrm{gen}}(\Gamma) ~ = ~ \mathrm{min}(|X|~|~X\in \mathrm{B}(\Gamma)) ~ \leq ~ \mathrm{sup}(|X|~|~X\in \mathrm{B}(\Gamma))
\end{equation} 
where the inequality is in fact an equality if $\Gamma$ admits no finite basis.  

Ranks different from the generating rank can also be considered. For instance, Buekenhout and Cohen \cite[Definition 5.3.1]{BC} propose the following: if $\Gamma$ admits a finite chain of subspaces of maximum length, that length is the rank of $\Gamma$; on the other hand, if no finite upper bound exists for the lengths of the chains of susbpaces of $\Gamma$, then take the symbol $\infty$ as the rank $\Gamma$ and say that $\Gamma$ has infinite rank. (Actually Buekenhout and Cohen define what they call the dimension; the rank is the dimension augmented by 1.) 

In their definition, Buekenhout and Cohen renounce to distinguish between different cases that can occur when the rank is infinite (as many authors do in cases like this). If we don't like this way of doing and want a sharper definition, then we can take the following as a rank of $\Gamma$:
\[\mathrm{rk}_{\mathrm{C}}(\Gamma) ~:= ~ \mathrm{sup}(\ell({\cal C})~|~ {\cal C}\in\mathfrak{C}(\Gamma))\] 
where $\mathfrak{C}(\Gamma)$ is the family of all chains of subspaces of $\Gamma$, namely sets of subspaces of $\Gamma$ totally ordered by inclusion, and $\ell({\cal C})$ is the length of a chain ${\cal C}\in\mathfrak{C}(\Gamma)$, namely $\ell({\cal C}) := |{\cal C}|-1$ where $|{\cal C}|$ is the cardinality of the set $\cal C$, with the usual convention that $|{\cal C}|-1 = |{\cal C}|$ when $|{\cal C}|$ is infinite. 

Suppose that $\Gamma$ admits a basis. As we shall prove later (Lemma \ref{lemma 2.1}), for every independent set $X$ there exists a chain ${\cal C}\in \mathfrak{C}(\Gamma)$ such that $|X| = \ell({\cal C})$. Hence 
\begin{equation}\label{intro 2}
\mathrm{sup}(|X|~|~X\in \mathrm{B}(\Gamma)) ~ \leq ~ \mathrm{sup}(|X|~|~X\in \mathrm{Ind}(\Gamma)) ~ \leq ~ \mathrm{rk}_{\mathrm{C}}(\Gamma).
\end{equation}
By (\ref{intro 1}) and (\ref{intro 2}) we immediately obtain the following inequality:
\begin{equation}\label{intro 3}
\mathrm{rk}_{\mathrm{gen}}(\Gamma)  ~ \leq ~ \mathrm{rk}_{\mathrm{C}}(\Gamma).
\end{equation}
In general (\ref{intro 3}) is a strict inequality, also because $\Gamma$ could admit bases of different cardinality or independent sets of cardinality  larger than $\mathrm{sup}(|X|~|~X\in \mathrm{B}(\Gamma))$ (see Subsection \ref{sec ex}, examples \ref{ex2} and \ref{ex4}). Oddities like these are possible because of the setting we have chosen, too weak for we can obtain anything sharper than inequalities. Its weakeness muddies the picture of the relations between $\mathrm{rk}_{\mathrm{gen}}(\Gamma)$ and $\mathrm{rk}_{\mathrm{C}}(\Gamma)$. In order to clear off that mud, henceforth we assume that $\Gamma$ satisfies the following property: 
\begin{itemize}
\item[(EP)] ({\bf Exchange Property}) {\em For $X\subseteq P$ and $x, y \in P$, if $y\in \langle X\cup\{x\}\rangle$ but $y\not\in \langle X\rangle$, then $x\in \langle X\cup\{y\}\rangle$.}
\end{itemize}
Then every $X\subseteq P$ contains a basis of its span $\langle X\rangle$. In particular, every $X \in \mathrm{Gen}(\Gamma)$ contains a minimal member of $\mathrm{Gen}(\Gamma)$ and the minimal members of $\mathrm{Gen}(\Gamma)$ are the same as the maximal members of $\mathrm{Ind}(\Gamma)$. Moreover, all bases have the same cardinality, necessarily equal to $\mathrm{rk}_{\mathrm{gen}}(\Gamma)$. We call $\mathrm{rk}_{\mathrm{gen}}(\Gamma)$ the {\em rank} of $\Gamma$ and we denote it by $\mathrm{rk}(\Gamma)$:
\[\mathrm{rk}(\Gamma) ~:= ~ \mathrm{rk}_{\mathrm{gen}}(\Gamma) ~ = ~ |X|,~~ \mbox{for any}~ X\in \mathrm{B}(\Gamma).\]
These claims can be proved by the same arguments commonly used to prove analogous claims in linear algebra. The next proposition is straightforward too (and well known): 

\begin{prop}\label{prop finite}
Assume {\rm (EP)} and let $\mathrm{rk}(\Gamma)$ be finite. Then $\mathrm{rk}_{\mathrm{C}}(\Gamma) = \mathrm{rk}(\Gamma)$.
\end{prop}

We might now be tempted to conjecture that the equality $\mathrm{rk}(\Gamma) = \mathrm{rk}_{\mathrm{C}}(\Gamma)$ holds in the infinite case as well, but this is false, as we shall prove in Section \ref{sec 2}. Explicitly:

\begin{theo}\label{th1}
Assuming {\rm (EP)}, let $\mathrm{rk}(\Gamma)$ be infinite and such that $2^\mathfrak{n} \leq \mathrm{rk}(\Gamma)$ for every cardinal number $\mathfrak{n} < \mathrm{rk}(\Gamma)$. Then 
\[\mathrm{rk}_{\mathrm{C}}(\Gamma) ~\geq ~ 2^{\mathrm{rk}(\Gamma)} ~ > ~ \mathrm{rk}(\Gamma).\] 
\end{theo} 
\begin{note}
{\rm The hypothesis that $2^\mathfrak{n}\leq \mathrm{rk}(\Gamma)$ whenever $\mathfrak{n} < \mathrm{rk}(\Gamma)$ is satisfied by any infinite $\mathrm{rk}(\Gamma)$ if we accept the Generalized Continuum Hypothesis $(\mathrm{GCH})$. Regardless of $(\mathrm{GCH})$, it holds if $\mathrm{rk}(\Gamma) = \aleph_0$.}
\end{note}

Theorem \ref{th1} makes it clear that, if we want to define the rank of $\Gamma$ by means of cardinalities of chains of subspaces, we cannot consider all possible such chains, namely $\mathrm{rk}_{\mathrm{C}}(\Gamma)$ is not the right candidate. The following is a way out: instead of considering the family $\mathfrak{C}(\Gamma)$ of all chains of subspaces of $\Gamma$, we can try the family $\mathfrak{W}(\Gamma)$ of all well ordered chains of subspaces of $\Gamma$, namely those chains where the inclusion relation defines a well ordering. Put 
\[\mathrm{rk}_{\mathrm{WO}}(\Gamma) ~ := ~ \mathrm{sup}(\ell({\cal C})~|~ {\cal C}\in \mathfrak{W}(\Gamma)).\]
Then the following holds, as we shall prove in Section \ref{sec 3}:

\begin{theo}\label{th2}
Assume {\rm (EP)}. Then $\mathrm{rk}(\Gamma) = \mathrm{rk}_{\mathrm{WO}}(\Gamma)$.
\end{theo}

So, $\mathfrak{W}(\Gamma)$ is the right family of chains to consider when we deal with ranks. Note that every finite totally ordered set is well ordered. When $\mathrm{rk}(\Gamma)$ is finite and (EP) holds then all chains of subspaces of $\Gamma$ are finite. In this case $\mathfrak{C}(\Gamma) = \mathfrak{W}(\Gamma)$. Proposition \ref{prop finite} is thus a special case of Theorem \ref{th2}. 

It is natural to ask what can be said on $\mathrm{rk}(\Gamma)$ and $\mathrm{rk}_{\mathrm{WO}}(\Gamma)$ when (EP) fails to hold. The next theorem, to be proved in Section \ref{sec 3 bis}, answers this question.  

\begin{theo}\label{th3}
In any case, $\mathrm{rk}(\Gamma) \leq  \mathrm{rk}_{\mathrm{WO}}(\Gamma)$.
\end{theo} 

Consequently, $\mathrm{rk}(\Gamma) \leq  \mathrm{rk}_{\mathrm{C}}(\Gamma)$ even if $\Gamma$ admits no bases.   

\begin{note}
{\rm Geometries where $\mathrm{rk}_{\mathrm{gen}}(\Gamma) < \mathrm{rk}_{\mathrm{WO}}(\Gamma)$ actually exist. See below, examples \ref{ex2} and \ref{ex4}.}
\end{note} 

\begin{note}\label{note EP}
{\rm Note that property (EP) deals with subspaces rather than lines. Geometries exist which have the same subspaces but different sets of lines. For instance, if a geometry $\Gamma$ admits lines of size $2$ and we remove some or all of them, or we take pairs of non-collinear points or proper subsets of lines of size at least $3$ as additional  lines, then we obtain a new geometry $\Gamma'$ with just the same subspaces as $\Gamma$ but a different set of lines. Property (EP) holds in $\Gamma$ if and only if it holds in $\Gamma'$.}
\end{note} 

\subsection{A few examples}\label{sec ex}

In the previous subsection, when disclaiming that certain properties hold in general, we promised suitable counterexamples. We now keep that promise.   

\begin{ex}\label{ex1}
{\rm Take set $\mathbb{N}$ of natural numbers as the set of points and the sets $L_n := \{kn~ |~ 0 \leq k \leq n\}$ as lines, for $n \geq 1$. Lines being defined in this way, two non-zero points $n$ and $m$ are collinear if and only if $n, m \leq d^2$, where $d$ is the greatest common divisor of $n$ and $m$. So, if $n$ and $m$ are collinear and $0 < m < n$, then $\langle m, n\rangle = \langle 0, n\rangle = \cup(L_u~|~ u~\mbox{divides}~n~\mbox{and}~ n\leq u^2\}$. 

Let $\Gamma$ be the geometry defined as above. It is not diffcult to see that a subset of $\mathbb{N}$ generates $\Gamma$ if and only if it contains a multiple of $n$ for every $n \neq 0$ and at least one pair of collinear points (which is certainly the case if it contains $0$). Let $X$ be such a set. Then $X$ always contains an element $m$ such that $X\setminus \{m\}$ still contains a pair of collinear points and multiples of every $n\neq 0$. So, $X\setminus \{m\}$ still generates $\Gamma$. Therefore no generating set of $\Gamma$ is minimal, namely $\Gamma$ admits no bases. Consequently, no maximal independent set can generate $\Gamma$. 

For instance, let $X_1$ be the set of prime numbers. Then $X_0 := X_1\cup\{0\}$ is a maximal independent set but it does not generate $\Gamma$. In fact $\langle X_0\rangle$ is the set of numbers of the form $n = pm$ for $p$ prime and $m\leq p$. If $n\in \mathbb{N}$ is not such a number then $n \not\in \langle X_0\rangle$ but $X_0\cup\{n\}$ is nevertheless a dependent set; indeed $p \in \langle (X_0\setminus \{p\})\cup\{n\}\rangle$ for at least one $p\in X_1$. 

Needless to say, (EP) fails to hold in $\Gamma$.}
\end{ex}  

\begin{ex}\label{ex2}
{\rm Let $\Gamma = (P,{\cal L})$ be defined as follows. For a point $a\in P$, the set $P\setminus \{a\}$ is partitioned in two subsets $B$ and $C$ of equal size $\mathfrak{n} = |B| = |C|$ and a bijection $f:B\rightarrow C$ is given. The set $B$ belongs to $\cal L$. The remaining lines of $\Gamma$ are the triples $\{a,b,f(b)\}$ for $b\in B$, the pairs $\{b,c\}$ with $b\in B$, $c\in C$ and $f(b) \neq c$ and all pairs of points of $C$. 

This geometry admits bases and all of its bases have size $3$. In fact $3 = \mathrm{rk}_{\mathrm{gen}}(\Gamma)$. On the other hand $C$ as well as the sets $C_b:= (C\setminus\{f(b)\})\cup\{b\}$ for $b\in B$ are maximal independent sets (whence $\mathfrak{n} = \mathrm{sup}(|X|~|~X\in \mathrm{Ind}(\Gamma))$), but none of them generates $\Gamma$. Accordingly, (EP) fails to hold in $\Gamma$.  

The sets $C$, all sets $C_b$ and their subsets are subspaces of $\Gamma$. Hence $\mathrm{rk}_{\mathrm{C}}(\Gamma) \geq \mathrm{rk}_{\mathrm{WO}}(\Gamma) \geq 1+\mathfrak{n}$. If $\mathfrak{n}$ is finite then $\mathrm{rk}_{\mathrm{C}}(\Gamma) = \mathrm{rk}_{\mathrm{WO}}(\Gamma) = 1+\mathfrak{n}$; when $\mathfrak{n}$ is infinite the inequality $\mathrm{rk}_{\mathrm{WO}}(\Gamma) \geq 1+\mathfrak{n}$ follows from Lemma \ref{lemma 2.1} of Section \ref{sec 3}. So, if $1+\mathfrak{n} > 3 = \mathrm{rk}_{\mathrm{gen}}(\Gamma)$, then $\mathrm{rk}_{\mathrm{C}}(\Gamma) \geq \mathrm{rk}_{\mathrm{WO}}(\Gamma) > \mathrm{rk}_{\mathrm{gen}}(\Gamma)$. 

Note that $\Gamma$ is a linear space \cite[Definition 2.5.13]{BC}. So, being a linear space is not enough for (EP) to hold. Conversely, (EP) is not sufficient for a geometry to be a linear space (see Remark \ref{note EP}).}
\end{ex} 

\begin{ex}\label{ex4}
{\rm Let $\Gamma$ be the polar space associated to a non-degenerate alternating form of $V(2n,\mathbb{F})$ for a field $\mathbb{F}$ and an integer $n > 1$ such that $2n < |\mathbb{F}|$. Then $\mathrm{rk}_{\mathrm{gen}}(\Gamma) = 2n$ but $\Gamma$ contains sets $X$ of pairwise non-collinear points of size $|X| = |\mathrm{PG}(1,\mathbb{F})| = 1+|\mathbb{F}|$. The hyperbolic lines of $\Gamma$ are sets like these, for instance. On the other hand, as no two points of $X$ are collinear, the set $X$ is independent. Moreover, $\mathrm{rk}_{\mathrm{WO}}(\Gamma) \geq |X|$ by Lemma \ref{lemma 2.1}. Hence $\mathrm{rk}_{\mathrm{WO}} \geq 1+|\mathbb{F}| > 2n = \mathrm{rk}_{\mathrm{gen}}(\Gamma)$.}
\end{ex}

\subsection{An outline of the rest of this paper}\label{outline}

Section 2 is devoted to the proofs of Theorems \ref{th1}, \ref{th2} and \ref{th3}. In Section 3 we discuss maximal chains of $\mathfrak{W}(\Gamma)$. Note that the union of a chain of well ordered chains is still a chain, but in general is not well ordered. So, we cannot ask Zorn's Lemma for help when looking for maximal members of $\mathfrak{W}(\Gamma)$. Nevertheless, as we shall show in Section 3, if $\Gamma$ satisfies (EP) then every well ordered chain of subspaces of $\Gamma$ is contained in a maximal well ordered chain. We shall also prove that a well ordered chain of subspaces of $\Gamma$ is maximal as a member of $\mathfrak{W}(\Gamma)$ only if it is maximal as a member of $\mathfrak{C}(\Gamma)$ and that all maximal well ordered chains of subspaces have the same length, equal to $\mathrm{rk}(\Gamma)$. 

In the last section of this paper (Section \ref{sec 5}) we shall turn back to the primary motivations of our investigation of chains of subspaces, namely the attempt to characterize the polar rank and the generating rank of a polar space by means of such chains. Explicitly, let $\Gamma$ be a non-degenerate polar space and define the polar rank $\mathrm{prk}(\Gamma)$ of $\Gamma$ as the least upper bound of the ranks of the singular subspaces of $\Gamma$, as usual. It is well known that if $\mathrm{prk}(\Gamma)$ is finite then $\mathrm{prk}(\Gamma)$ is equal to the maximal length of a chain of singular subspaces. On the other hand, when  $\mathrm{prk}(\Gamma)$ is infinite it can happen that chains of singular subspaces exist of length greater than $\mathrm{prk}(\Gamma)$. This is not surprising, in view of Theorem \ref{th1}. Theorem \ref{th2} suggests a way out:  $\mathrm{prk}(\Gamma)$ is equal to the least upper bound of the well ordered chains of singular subspaces, as it follows from Theorem \ref{th2} applied to each of the maximal singular subspaces of $\Gamma$. 

Turning to the generating rank, let $\Gamma$ be a non-degenerate polar space of polar rank at least 2. Property (EP) fails to hold in $\Gamma$. Hence we cannot apply Theorem \ref{th2}. In fact, in general $\mathrm{rk}_{\mathrm{gen}}(\Gamma) < \mathrm{rk}_{\mathrm{WO}}(\Gamma)$ (compare Example \ref{ex4}). However, as shown in \cite{CGP}, we can bypass this obstacle by considering only subspaces which contain a pair of mutually disjoint maximal singular subspaces. As in \cite{CGP}, we call these subspaces {\em nice subspaces}. Let $\mathfrak{W}_{\mathrm{nice}}(\Gamma)$ be the family of well ordered chains of nice subspaces and assume that $\Gamma$ has finite polar rank $n = \mathrm{prk}(\Gamma)$. Then the following holds for a large class of polar spaces, which includes classical polar spaces defined over a commutative division ring or a division ring of characteristic different from 2: 
\[\mathrm{rk}_{\mathrm{gen}}(\Gamma) ~ = ~ 2n + \mathrm{sup}(\ell({\cal C})~|~ {\cal C}\in\mathfrak{W}_{\mathrm{nice}}(\Gamma)).\]
        
\section{Proof of Theorems \ref{th1}, \ref{th2} and \ref{th3}}

\subsection{Proof of Theorem \ref{th1}}\label{sec 2}

\begin{lemma}\label{lemma 1.1}
Let $U$ be an infinite set such that $2^\mathfrak{n}\leq |U|$ for every cardinal number $\mathfrak{n} < |U|$. Then there exists a chain $\cal C$ of subsets of $U$ such that $|{\cal C}| = 2^{|U|}$.
\end{lemma}
{\bf Proof.} Let $\omega$ be the least ordinal number such that $|\omega| = |U|$. Regarded $\omega$ as the same as the well ordered set of all ordinal numbers $\delta < \omega$, we can consider the set $\{0,1\}^\omega$ of all mappings from $\omega$ to $\{0,1\}$, namely transfinite sequences of length $\omega$ with $0$ or $1$ as entries. We can order $\{0,1\}^\omega$ lexicographically, by declaring that $f < g$ for two mappings $f, g \in \{0,1\}^\omega$ if and only if there exists an ordinal number $\gamma < \omega$ such that $f(\gamma) = 0$, $g(\gamma) = 1$ and $f(\delta) = g(\delta)$ for every $\delta < \gamma$. In this way a total order is defined on $\{0,1\}^\omega$. 

For $\gamma < \omega$ let $\{0,1\}^{\omega|\gamma}$ be the set of all mappings $f\in\{0,1\}^\omega$ such that $f(\delta) = 0$ for any $\gamma \leq \delta < \omega$. Clearly $|\{0,1\}^{\omega|\gamma}| = |\{0,1\}^\gamma|$. Moreover $2^{|\gamma|} \leq |\omega| = |U|$ by the hypotheses assumed on $|U|$ and since $|\gamma| < |\omega|$ by the choice of $\omega$. It follows that the set $\{0,1\}^{\omega|*}  := \cup_{\gamma < \omega}\{0,1\}^{\omega|\gamma} \subset \{0,1\}^\omega$ 
has the same cardinality as $\omega$, hence the same as $U$. Accordingly, there exists a bijection $\chi:U\rightarrow \{0,1\}^{\omega|*} \subset \{0,1\}^\omega$. 

For every $f\in \{0,1\}^\omega$ let $U_f := \{x\in U~|~ \chi(x) < f\}$. It is easily seen that for $f, g \in \{0,1\}^\omega$ we have  $U_f\subset U_g$ if and only if $f < g$ in the lexicographic order previously defined on $\{0,1\}^\omega$. Hence the set ${\cal C} :=  \{U_f~|~ f \in \{0,1\}^\omega\}$ 
is a chain of subsets of $U$ of cardinality $|{\cal C}| = |\{0,1\}^\omega| = 2^{|U|}$. \hfill $\Box$

\begin{note}
{\rm The chain $\cal C$ constructed in the proof of Lemma \ref{lemma 1.1} is not maximal. Indeed the sets $\overline{U}_f  = U_f\cup\{\chi^{-1}(f)\}$ for $f\in \{0,1\}^{\omega|*}$ are missing in $\cal C$. If we insert them too in $\cal C$ then we obtain a maximal chain, say $\overline{\cal C}$. Clearly, $\ell(\overline{\cal C}) = \ell({\cal C})$.}
\end{note}

Let $\Gamma$ be a point-line geometry. Without assuming (EP), suppose that $\Gamma$ admits a basis and let $U\in B(\Gamma)$ be a basis of $\Gamma$.

\begin{prop}\label{prop 1.2}
Suppose that $U$ is infinite and $2^\mathfrak{n} \leq |U|$ for every cardinal number $\mathfrak{n} < |U|$. Then $\mathrm{rk}_{\mathrm{C}}(\Gamma) \geq 2^{|U|}$. 
\end{prop}
{\bf Proof.} Since $U$ is a basis, the span operator $\langle .\rangle$ induces an inclusion preserving injective mapping from the boolean lattice of subsets of $U$ to the subspace lattice of $\Gamma$. In particular, it maps isomorphically every chain of subsets of $U$ onto a chain of subspaces of $\Gamma$. The conclusion follows from Lemma \ref{lemma 1.1}. \hfill $\Box$ 
\\

Theorem \ref{th1} immediately follows from Proposition \ref{prop 1.2}.   
 
\subsection{Proof of Theorem \ref{th2}}\label{sec 3}

Let $\Gamma = (P, {\cal L})$ be a point-line geometry, possibly not satisfying the Exchange Property (EP). Note that $\mathrm{Ind}(\Gamma) \neq \emptyset$ even if $\Gamma$ admits no basis. Indeed $\emptyset$, the singletons of the points and the pairs of distinct points of $\Gamma$ are independent sets. By Zorn's Lemma, every independent set is contained in a maximal independent set. However, if (EP) fails to hold in $\Gamma$ then a maximal independent set need not generate $\Gamma$. In particular, $\Gamma$ could admit no bases. Nevertheless, the following number is well defined:
\[\mathrm{rk}_{\mathrm{ind}}(\Gamma) ~ := ~ \mathrm{sup}(|X|~|~ X\in \mathrm{Ind}(\Gamma)). \]

\begin{lemma}\label{lemma 2.1}
We have $\mathrm{rk}_{\mathrm{ind}}(\Gamma) \leq \mathrm{rk}_{\mathrm{WO}}(\Gamma)$. 
\end{lemma}
{\bf Proof.} Given $X\in \mathrm{Ind}(\Gamma)$, let $\xi = (x_\delta)_{\delta < \omega}$ be a well ordering of $X$, namely a bijective mapping as follows:
\[\begin{array}{rccl}
\xi : & \omega = \{\delta\}_{\delta < \omega}  & \stackrel{1-1}{\longrightarrow} & X \\
 & \delta < \omega & \rightarrow & x_\delta = \xi(\delta)
\end{array}\]
For every $\gamma \leq \omega$ put $S_\gamma := \langle x_\delta\rangle_{\delta < \gamma}$. As $X$ is independent, if $0 \leq \gamma_1 < \gamma_2 \leq \omega$ then $S_{\gamma_1} \subset S_{\gamma_2}$. Therefore ${\cal C}_{X,\xi} := \{S_\gamma\}_{\gamma \leq \omega}$ is a well ordered chain of length $\ell({\cal C}_{X,\xi}) = |\omega| = |X|$. \hfill $\Box$

\bigskip

Conversely, let ${\cal C} = \{S_\delta\}_{\delta \leq \omega}$ be a well ordered chain of subspaces of $\Gamma$ with $S_0 = \emptyset$. As usual when dealing with well ordered sets, we assume that $S_{\gamma_1} \subset S_{\gamma_2}$ if and only if $\gamma_1 < \gamma_2$. For every $\delta < \omega$ we choose a point $x_\delta\in S_{\delta+1}\setminus S_\delta$ and we put $X_{{\cal C},\xi} := \{x_\delta\}_{\delta < \omega}$, where 
\[\begin{array}{rccl}
\xi: & \omega  & \longrightarrow & \cup_{\gamma \leq \omega}S_\gamma \\
 & \delta < \omega & \rightarrow & x_\delta \in S_{\delta+1}\setminus S_\delta
\end{array}\]
is the function implicitly defined by those choices.  

\begin{lemma}\label{lemma 2.2}
Assume {\rm (EP)}. Then the set $X_{{\cal C},\xi}$ is independent and $|X_{{\cal C},\xi}| = \ell({\cal C})$. 
\end{lemma}
{\bf Proof.} The equality $|X_{{\cal C},\xi}| = \ell({\cal C})$ is obvious. By contradiction, suppose that $X_{{\cal C},\xi}$ is not an independent set. Let $\gamma\leq \omega$ be the smallest ordinal number such that $\{x_\delta\}_{\delta < \gamma}$ is dependent and let $\delta_0 < \gamma$ be such that $x_{\delta_0} \in \langle x_\delta ~ |~ \delta_0 \neq \delta < \gamma\rangle$. Then $x_{\delta_0}\in \langle x_{\delta_i}\rangle_{i=1}^n$ for a finite subset $\{\delta_1,..., \delta_n\}\subseteq \{\delta < \gamma~|~ \delta \neq \delta_0\}$. We can assume that $\delta_1 < \delta_2 < ... < \delta_n$. Clearly $\eta := \mathrm{max}(\delta_0, \delta_n) < \gamma$. Therefore $\gamma = \eta+1$ by the minimality of $\gamma$. 

Suppose firstly that $\eta = \delta_0 > \delta_n$. Then $\langle x_{\delta_1},...,x_{\delta_n}\rangle\subseteq S_{\delta_n+1}\subset S_{\delta_0+1}$. However $x_{\delta_0} \in S_{\delta_0+1}\setminus S_{\delta_0} \subseteq S_{\delta_0+1}\setminus S_{\delta_n+1}$. This contradicts the hypothesis that $x_{\delta_0} \in \langle x_{\delta_1},...,x_{\delta_n}\rangle$. Therefore $\delta_0 < \delta_n = \eta$. 

Neither $x_{\delta_0}$ nor $x_{\delta_n}$ belong to $S := \langle x_{\delta_1},....,x_{\delta_{n-1}}\rangle$. Indeed $x_{\delta_0}\not\in S$ by the minimality of $\gamma = \delta_n+1$ and $x_{\delta_n}\not \in S$ since $S\subseteq S_{\delta_n}$ and $x_{\delta_n}\in S_{\delta_n+1}\setminus S_{\delta_n}$ by definition. On the other hand, $x_{\delta_0}\in \langle S\cup\{x_{\delta_n}\}\rangle$ by assumption. Hence (EP) forces $x_{\delta_n}\in\langle S\cup \{x_{\delta_0}\}\rangle \subseteq S_{\delta_n}$ while $x_{\delta_n}\in S_{\delta_n+1}\setminus S_{\delta_n}$ by definition. We have reached a final contradiction.  \hfill $\Box$ 

\bigskip  

By Lemmas \ref{lemma 2.1} and \ref{lemma 2.2} we immediately obtain the following:

\begin{prop}\label{prop 2.3}
If {\rm (EP)} holds in $\Gamma$ then $\mathrm{rk}_{\mathrm{ind}}(\Gamma) = \mathrm{rk}_{\mathrm{WO}}(\Gamma)$.
\end{prop}

However, if (EP) holds then $\mathrm{rk}_{\mathrm{ind}}(\Gamma) = \mathrm{rk}(\Gamma)$. Hence $\mathrm{rk}_{\mathrm{WO}}(\Gamma) = \mathrm{rk}(\Gamma)$, as claimed in Theorem \ref{th2}.  

\subsection{Proof of Theorem \ref{th3}}\label{sec 3 bis}

Let $X\in \mathrm{Gen}(\Gamma)$. We can give $X$ a well ordering $X = \{x_\delta\}_{\delta < \omega}$. For every $\gamma \leq \omega$, put $S_\gamma := \langle x_\delta\rangle_{\delta < \gamma}$. The set ${\cal C} := \{S_\gamma\}_{\gamma \leq \omega}$ is a chain of subspaces with $S_0 = \emptyset$ and $S_\omega = \langle X\rangle = P$. Clearly, if $\gamma_1 < \gamma_2$ then $S_{\gamma_1}\subseteq S_{\gamma_2}$ but, as $X$ is not assumed to be independent, it could happen that $S_{\gamma_1} = S_{\gamma_2}$. Nevertheless:

\begin{lemma}\label{lemma 2.4}
The chain $\cal C$ is well ordered.
\end{lemma}
{\bf Proof.} Let $\{S_\gamma\}_{\gamma\in U}\subseteq{\cal C}$, for $U \subseteq \omega+1 = \{\gamma\}_{\gamma \leq \omega}$.
Let $\gamma_0$ be the minimum of $U$. Then $S_{\gamma_0}$ is the minimum of $\{S_\gamma\}_{\gamma \in U}$. \hfill $\Box$

\bigskip

By Lemma \ref{lemma 2.4}, there exists an ordinal number $\omega' \leq \omega$ such that ${\cal C} = \{S'_\eta\}_{\eta \leq \omega'}$, where the indices $\eta \leq \omega'$ are chosen in such a way that
$S'_{\eta_1} \subset S'_{\eta_2}$ if and only if $\eta_1 < \eta_2$. For every $\eta\leq \omega'$ let $f(\eta)$ be the least $\gamma\leq \omega$ such that $S_\gamma = S'_{\eta}$ and $g(\eta)$ the least $\delta < \omega$ such that $x_{g(\eta)}\not \in S'_\eta$. Then $f$ and $g$ are injective morphisms of posets from $\omega'+1 = \{\eta\}_{\eta \leq \omega'}$ to $\omega+1$ and $f(\eta) \leq g(\eta)$ for every $\eta \leq \omega'$. Clearly ${\cal C} = \{S_{g(\eta)}\}_{\eta \leq \omega'}$. In particular  $S'_{\omega'} = S_{f(\omega')} = S_\omega$, but possibly $f(\omega') < \omega$.     

\begin{lemma}\label{lemma 2.5}
We have $S'_{\eta} = \langle x_{g(\varepsilon)}\rangle_{\varepsilon < \eta}$ for every $\eta\leq \omega'$. 
\end{lemma}
{\bf Proof.} The proof is by induction. As $S'_0 = S_0 = \emptyset$, if $\eta = 0$ there is nothing to prove. Let $\eta > 0$ and $\gamma_1 = f(\eta)$. Suppose firstly that $\gamma_1$ is a limit ordinal. Then, by definition of $f$ and $g$ and the inductive hypothesis, we have  
\[S_{\gamma_1} = \cup_{\gamma < \gamma_1}S_\gamma = \cup_{\varepsilon < \eta}S_{f(\varepsilon)} = \cup_{\varepsilon < \eta}S'_\varepsilon = \cup_{\varepsilon < \eta}\langle x_{g(\zeta)}\rangle_{\zeta < \varepsilon} = \langle x_{g(\varepsilon)}\rangle_{\varepsilon < \eta}.\]
Hence $S'_\eta = S_{\gamma_1} =  \langle x_{g(\varepsilon)}\rangle_{\varepsilon < \eta}$, as claimed. On the other hand, let $\gamma_1 = \gamma_0+1$.  Then $S_{\gamma_0} \subset S'_\eta = S_{\gamma_1}$, since $\gamma_1 = f(\eta)$. Also, if $\eta' < \eta$ is such that $S'_{\eta'} = S_{\gamma_0}$, then $f(\eta') \leq \gamma_0$ but $\eta = \eta'+1$. Therefore $S_{\gamma_0} = S'_{\eta'} = \langle x_{g(\varepsilon)}\rangle_{\varepsilon < \eta'}$ by the iductive hypothesis, while $S_{\gamma_1} = \langle S_{\gamma_0}\cup\{x_{\gamma_0}\}\rangle$. It follows that $g(\gamma_0) = \gamma_0$. Hence $S'_\eta = \langle x_{g(\varepsilon)}\rangle_{\varepsilon < \eta}$ in this case too. \hfill $\Box$

\bigskip

Put $X' = \{x_{\delta(\gamma)}\}_{\gamma \leq \omega'} \subseteq X$. By Lemma \ref{lemma 2.5} and since $S'_{\omega'} = S_\omega = P$, the set $X'$ generates $\Gamma$. Hence $|X'| \geq \mathrm{rk}_{\mathrm{gen}}(\Gamma)$. On the other hand, $|X'| = |\omega'| = |{\cal C}| \leq \mathrm{rk}_{\mathrm{WO}}(\Gamma)$. Hence $\mathrm{rk}_{\mathrm{gen}}(\Gamma) \leq \mathrm{rk}_{\mathrm{WO}}(\Gamma)$, as claimed in Theorem \ref{th3}.  

\section{Maximal well ordered chains}\label{sec 4}

Throughout this section $\Gamma = (P, {\cal L})$ is a point-line geometry satisfying (EP). Given $X\in \mathrm{Ind}(\Gamma)$ and a well ordering $\xi = (x_\delta)_{\delta < \omega}$ of $X$, let ${\cal C}_{X,\xi} = \{S_\gamma\}_{\gamma \leq \omega}$ be the well ordered chain defined as in the proof of Lemma \ref{lemma 2.1}. Note that $S_0 = \emptyset$. If $\gamma$ is a limit ordinal then $S_\gamma := \cup_{\delta < \gamma}S_\delta$ and if $\gamma = \delta+1$ then $S_\gamma := \langle S_\delta\cup\{x_\delta\}\rangle$.

\begin{lemma}\label{lemma 3.1}
The following are equivalent:

\begin{itemize}
\item[$(1)$] The set $X$ is a basis of $\Gamma$;
\item[$(2)$] $S_\omega = P$;
\item[$(3)$] The chain ${\cal C}_{X,\xi}$ is maximal as a member of $\mathfrak{C}(\Gamma)$; 
\item[$(4)$] The chain ${\cal C}_{X,\xi}$ is maximal as a member of $\mathfrak{W}(\Gamma)$.
\end{itemize}
\end{lemma}
{\bf Proof.} By definition, $S_\gamma = \langle x_\delta\rangle_{\delta < \gamma}$. In particular, $S_\omega = \langle X\rangle$. The equivalence of (1) and (2) is obvious. Clearly, (3) implies (4) and (4) in turn implies (2). Indeed, if $S_\omega\subset P$ then ${\cal C}_{X,\xi}\cup\{P\}$ is a well ordered chain properly containing ${\cal C}_{X,\xi}$; thus ${\cal C}_{X,\xi}$ cannot be maximal in $\mathfrak{W}(\Gamma)$. 

The implication $(1)\Rightarrow(3)$ remains to be proved. Suppose that ${\cal C}_{X,\xi}$ is not a maximal chain. Then there exists a subspace $S$ such that for every $\gamma \leq \omega$ either $S_\gamma \subset S$ or $S \subset S_\gamma$. Let $J := \{\gamma \leq \omega~|~ S \subset S_\gamma\}$.  
Suppose firstly that $J \neq \emptyset$ and let $\gamma$ be the minimum element of $J$. So, $S_\gamma \supset S \supset S_\delta$ for every $\delta < \gamma$. If $\gamma$ is a limit ordinal this cannot be, since in this case $S_\gamma = \cup_{\delta < \gamma}S_\delta$. Therefore $\gamma$ admits a predecessor, say $\gamma = \delta+1$. Then $S_\delta \subset S \subset S_{\delta+1} = S_\gamma$. However $S_\gamma = \langle S_\delta\cup \{x_\delta\}\rangle$. Hence, by (EP), no subspace exists which properly contains $S_\delta$ and is properly contained in $\langle S_\delta\cup \{x_\delta\}\rangle$. We have reached a contradiction. Therefore $J = \emptyset$. Accordingly, $S_\omega \subset S$. Hence $S_\omega \neq P$ and $X$ is not a basis.  \hfill $\Box$ 

\bigskip

The following is implicit in the proof of Lemma \ref{lemma 3.1}:

\begin{cor}\label{cor 3.2}
If $\cal C$ is a chain containing ${\cal C}_{X,\xi}$ then all members of ${\cal C}\setminus {\cal C}_{X,\xi}$ properly contain the largest element $S_\omega$ of ${\cal C}_{X,\xi}$.
\end{cor}

Conversely, given a well ordered chain ${\cal C} = \{S_\gamma\}_{\gamma \leq \omega}$ with $S_0 = \emptyset$, let $X_{{\cal C},\xi}$ be constructed as in the second part of Section \ref{sec 3}. As $\Gamma$ is assumed to satisfy (EP), the set $X_{{\cal C},\xi}$ is independent by Lemma \ref{lemma 2.2}. 

\begin{theo}\label{prop 3.3}
The following are equivalent:
\begin{itemize} 
\item[$(1)$] The chain $\cal C$ is maximal as a member of $\mathfrak{C}(\Gamma)$; 
\item[$(2)$] The chain $\cal C$ is maximal as a member of $\mathfrak{W}(\Gamma)$;
\item[$(3)$] The set $X_{{\cal C},\xi}$ is a basis of $\Gamma$. 
\end{itemize}
\end{theo}
{\bf Proof.} Recall that $X_{{\cal C},\xi} = \{x_\delta\}_{\delta < \omega}$ where $x_\delta = \xi(\delta) \in S_{\delta+1}\setminus S_{\delta}$. For every $\gamma \leq \omega$ put $S_{\gamma,\xi} := \langle x_\delta\rangle_{\delta < \gamma}$. So, $S_{\gamma, \xi} \subseteq S_\gamma$ and the chain ${\cal C}_\xi := \{S_{\gamma, \xi}\}_{\gamma \leq \omega}$ is the same as the chain ${\cal C}_{X,\xi}$ constructed from $X = X_{{\cal C},\xi}$ by exploiting the same function $\xi$ used to construct $X$ from $\cal C$.    

Let $J := \{\gamma \leq \omega ~|~ S_{\gamma,\xi} \subset S_\gamma\}$. If $J = \emptyset$ then ${\cal C}_\xi = {\cal C}$. In this case the equivalence of (1), (2) and (3) follows from Lemma \ref{lemma 3.1}. Suppose $J \neq \emptyset$ and let $\gamma$ be the minimum of $J$. Note that $0 < \gamma$ since $S_0 = \emptyset$ by assumption and $S_{0, \xi} = \emptyset$ by definition. We have $S_{\delta, \xi} = S_\delta$ for every $\delta < \gamma$. If $\gamma$ is a limit ordinal then $S_{\gamma,\xi} = \cup_{\delta < \gamma}S_{\delta,\xi} = \cup_{\delta < \gamma}S_\delta = S_\gamma$, contrary to the hypothesis that $S_{\gamma, \xi} \subset S_\gamma$. Therefore $\{\delta~|~\delta < \gamma\}$ admits a maximum $\eta$ and $\gamma = \eta+1$. Hence $S_\eta = S_{\eta,\xi} \subset S_{\gamma,\xi} \subset S_\gamma$. 

We can insert $S_{\gamma, \xi}$ in $\cal C$ between $S_\eta$ and $S_\gamma$, thus obtaining a well ordered chain ${\cal C}' = (S'_\delta)_\delta\leq \omega'$ properly containing $\cal C$. Explicitly, if at least one limit ordinal occurs between $\gamma$ and $\omega$ then $\omega' = \omega$, otherwise $\omega' = \omega+1$. As for the subspaces $S'_\delta$, they are defined as follows. If $\delta < \gamma$ or $\delta \leq \omega$ is such that $\delta > \gamma+n$ for every finite ordinal $n$, then we put $S'_\delta := S_\delta$. (Note that an ordinal $\delta \leq \omega$ such that $\delta > \gamma+n$ for every $n$ exists precisely when $\omega' = \omega$.) Moreover $S'_\gamma := S_{\gamma, \xi}$ and $S'_{\gamma+n+1} := S_{\gamma+n}$ for every finite ordinal $n$. 

The chain $\cal C$, being contained in ${\cal C}'$ is not maximal in $\mathfrak{W}(\Gamma)$; even more so, it cannot be maximal in $\mathfrak{C}(\Gamma)$. Moreover, we can extend the mapping $\xi$ to a mapping $\chi:\omega'\rightarrow \cup_{\delta\leq\omega'}S'_\delta$ as follows. If $\delta < \gamma$ of $\delta > \gamma+n$ for every finite ordinal $n$ then we set $\chi(\delta) := \xi(\delta) = x_\delta$. We put $\chi(\gamma+n+1) = \xi(\gamma+n) = x_{\gamma+n}$ for every finite ordinal $n$ and we choose an element of $S'_{\gamma+1}\setminus S'_\gamma = S_\gamma\setminus S_{\gamma,\xi}$ as $\chi(\gamma)$. The set $X_{{\cal C}',\chi} := \{\chi(\delta)\}_{\delta\leq \omega'}$ is independent by Lemma \ref{lemma 2.2} and properly contains $X_{{\cal C},\xi}$. Consequently, $X_{{\cal C},\xi}$ is not a basis of $\Gamma$.

Summarizing, when $J = \emptyset$ then (1), (2) and (3) are equivalent by Lemma \ref{lemma 3.1}. When $J\neq \emptyset$ then all claims (1), (2) and (3) are false, whence trivially equivalent. So, (1), (2) and (3) are equivalent in any case. \hfill $\Box$ 

\bigskip

The next corollary immediately follows from Theorem \ref{prop 3.3}. 

\begin{cor}\label{cor 3.4}
If a well ordered chain of subspaces of $\Gamma$ is maximal as a member of $\mathfrak{W}(\Gamma)$ then it is maximal in $\mathfrak{C}(\Gamma)$ too.
\end{cor} 

\begin{cor}\label{cor 3.5} 
A well ordered chain ${\cal C} = \{S_\gamma\}_{\gamma\leq \omega}$ is maximal if and only if all of the following hold: 

\medskip

$(1)$ ~ $S_\gamma$ is a maximal subspace of $S_{\gamma+1}$ for every $\gamma < \omega$;

$(2)$ ~ $\cup_{\delta < \gamma}S_\delta = S_\gamma$ for every limit ordinal $\gamma \leq \omega$;

$(3)$ ~ $S_0 = \emptyset$;

$(4)$ ~ $S_\omega = P$. 
\end{cor}
{\bf Proof.}  With $J$ as in the proof of Theorem \ref{prop 3.3}, conditions (1), (2) and (3) characterize the case $J = \emptyset$ and, in this case, ${\cal C} = {\cal C}_{X,\xi}$. By Lemma \ref{lemma 3.1}, the chain ${\cal C}_{X,\xi}$ ($= {\cal C}$) is maximal if and only if (4) holds in it. \hfill $\Box$ 

\begin{cor}\label{cor 3.6}
All maximal well ordered chains of subspaces of $\Gamma$ have length equal to $\mathrm{rk}(\Gamma)$.
\end{cor}
{\bf Proof.} As shown in the proof of Corollary \ref{cor 3.5}, if ${\cal C}\in \mathfrak{W}(\Gamma)$ is maximal then ${\cal C} = {\cal C}_{X,\xi}$ and $X$ is a basis of $\Gamma$. However $\ell({\cal C}_{X,\xi}) = |X|$ and $|X| = \mathrm{rk}(\Gamma)$ because $X\in \mathrm{B}(\Gamma)$ and (EP) holds by assumption. Hence $\ell({\cal C}) = \mathrm{rk}(\Gamma)$. \hfill $\Box$    

\begin{theo}\label{theo 3.7}
Every well ordered chain of subspaces of $\Gamma$ is contained in a maximal well ordered chain.
\end{theo}
{\bf Proof.} Recall that $\Gamma$ satsifies (EP), by assumption. Hence the following holds:  

\begin{itemize} 
\item[$(\ast)$] every independent subset of a subspace $S$ is contained in a basis of $S$. 
\end{itemize}
Let ${\cal C} \in \mathfrak{W}(\Gamma)$. As we want to prove that $\cal C$ is contained in a maximal member of $\mathfrak{W}(\Gamma)$, there is loss in assuming that $\cal C$ contains $\emptyset$ and the full point-set $P$ of $\Gamma$. Indeed every maximal well ordered chain contains these two subspaces. So, ${\cal C} = \{S_\gamma\}_{\gamma\leq \omega}$ with $S_0 = \emptyset$ and $S_\omega = P$.

We shall now define a chain $\{X_\gamma\}_{\gamma\leq\omega}$ of independent sets such that if $\delta < \gamma \leq \omega$ then $X_\delta \subset X_\gamma$ and $X_\gamma$ is a basis of $S_\gamma$, for every $\gamma \leq \omega$. We put $X_0 = \emptyset$ and we go on by induction. Let $\gamma > 0$ and assume to have already defined a basis $X_\delta$ of $S_\delta$ for every $\delta < \gamma$ in such a way that if $\delta < \eta < \gamma$ then $X_\delta \subset X_\eta$. The union $X'_\gamma := \cup_{\delta < \gamma}X_\delta$ is an independent subset of $S_\gamma$. By $(*)$, the set $X'_\gamma$ is contained in a basis $X_\gamma$ of $S_\gamma$. 

Clearly, $X_\omega$ is a basis of $\Gamma$. We can also give each of the sets $X_\gamma$ a well ordering $\xi_\gamma$ in such a way that if $\delta < \gamma$ then $\xi_\gamma$ induces $\xi_\delta$ on $X_\delta$. This too can be done by induction. As $X_0 = \emptyset$ take the empty order as $\xi_0$. Assume to have defined $\xi_\delta$ for every $\delta < \gamma$ in such a way that $\xi_\eta$ induces $\xi_\delta$ on $X_\delta$ if $\delta < \eta < \gamma$. Let $\chi_\gamma$ be the well ordering thus defined on $X'_\gamma$. Explicitly, if $\gamma$ is a limit ordinal then $\chi_\gamma$ is the limit $\lim_{\delta < \gamma}\xi_\delta$ of the sequence of well orders $(\xi_\delta)_{\delta < \gamma}$. On the other hand, if $\gamma = \eta+1$ then $\chi_\gamma = \xi_\eta$. 

Put $Y_\gamma := X_\gamma\setminus X'_\gamma$. Note that when $\gamma$ is a limit ordinal it can happen that $X'_\gamma$ is a basis of $S_\gamma$. If this is the case then $X_\gamma = X'_\gamma$ and $Y_\gamma = \emptyset$. Choose a well ordering $\zeta_\gamma$ of $Y_\gamma$ and define $\xi_\gamma$ as the sum $\xi_\gamma = \chi_\gamma + \zeta_\gamma$ of $\chi_\gamma$ and $\zeta_\gamma$. Explicitly, $\xi_\gamma$ induces $\chi_\gamma$ and $\zeta_\gamma$ on $X'_\gamma$ and $Y_\gamma$ respectively and every element of $Y_\gamma$ follows all elements of $X'_\gamma$ in the order $\xi_\gamma$.

The basis $X_\omega$ of $\Gamma$ thus gets a well ordering $\xi_\omega: \omega'\rightarrow X$ which induces $\xi_\gamma$ on $X_\gamma$ for every $\gamma \leq \omega$. Clearly $\omega' \geq \omega$. Let ${\cal C}_{X_\omega,\xi_\omega}$ be the well ordered chain constructed as in the proof of Lemma \ref{lemma 2.1} for $(X,\xi) = (X_\omega,\xi_\omega)$. The chain ${\cal C}_{X_\omega,\xi_\omega}$ is maximal by Lemma \ref{lemma 3.1} and contains $\cal C$, by construction.  \hfill $\Box$ 

\bigskip

So far for well-ordered chains. When $\Gamma$ is finitely generated then all chains of subspaces of $\Gamma$ are finite, hence well-ordered. On the other hand, when $\mathrm{rk}(\Gamma)$ is infinite non-well ordered chains always exist, even maximal ones. By Corollary \ref{cor 3.4}, a non-well ordered maximal member of $\mathfrak{C}(\Gamma)$ contains no maximal members of $\mathfrak{W}(\Gamma)$. So, Corollary \ref{cor 3.6} is of no use to draw any conclusion on the possible length of a non-well ordered maximal member of $\mathfrak{C}(\Gamma)$. We conjecture the following:  

\begin{conj}\label{conj1}
We have $\mathrm{rk}(\Gamma) \leq \ell({\cal C}) \leq 2^{\mathrm{rk}(\Gamma)}$ for every maximal ${\cal C}\in \mathfrak{C}(\Gamma)$.
\end{conj}

Recall that in Theorem \ref{th1} we assume that $2^\mathfrak{n}\leq \mathrm{rk}(\Gamma)$ for every $\mathfrak{n} < \mathrm{rk}(\Gamma)$.

\begin{conj}
When $\mathrm{rk}(\Gamma)$ is infinite, chains of subspaces of length at least $2^{\mathrm{rk}(\Gamma)}$ (exactly $2^{\mathrm{rk}(\Gamma)}$ if Conjecture {\rm \ref{conj1}} holds true) exist even if $2^\mathfrak{n} > \mathrm{rk}(\Gamma)$ for some $\mathfrak{n} < \mathrm{rk}(\Gamma)$.
\end{conj} 

\section{Two problems in the theory of polar spaces}\label{sec 5} 

\subsection{The polar rank of a polar space} 

Let $\Gamma$ be a non-degenerate polar space. Its singular subspaces are projective spaces (Buekenhout and Shult \cite{BS}, also Buekenhout and Cohen \cite{BC} or Shult \cite{Shult}) and each of them is contained in a maximal singular subspace. When at least one of the maximal singular subspaces of $\Gamma$ has finite rank (the rank $\mathrm{rk}(S)$ of a projective space $S$ being its generating rank), then all of them have the same rank (see e.g. \cite{BC}, \cite{BS} or \cite{Shult}); that common rank is usually called the rank of $\Gamma$. However, following Shult \cite{Shult}, we prefer to call it the {\em polar rank}, in order to avoid any confusion with the generating rank $\mathrm{rk}_{\mathrm{gen}}(\Gamma)$ of $\Gamma$ as a point-line geometry, which is larger than the polar rank. Henceforth we denote the polar rank of $\Gamma$ by the symbol $\mathrm{prk}(\Gamma)$. So,
\begin{equation}\label{def1}
\mathrm{prk}(\Gamma) ~ = ~ \mathrm{sup}(\mathrm{rk}(S)~|~ S~\in \mathrm{Sing}(\Gamma)).
\end{equation}
where $\mathrm{Sing}(\Gamma)$ stands for the collection of all singular subspaces of $\Gamma$. Nearly all authors define $\mathrm{prk}(\Gamma)$ according to (\ref{def1}), but only when all singular subspaces of $\Gamma$ have finite rank. When $\Gamma$ also admits singular subspaces of infinite rank they cut short by stating that $\mathrm{prk}(\Gamma) = \infty$. (Compare the definition of the dimension of a line-space by Buekenhout and Cohen \cite[Definition 5.3.1]{BC}, discussed in Section \ref{sec 1}.) Let's call this convention the $\infty$-{\em convention}. 

However equality (\ref{def1}) also makes sense when $\Gamma$ admits singular subspaces of infinite rank. Since I don't like the $\infty$-convention so much, I prefer to take (\ref{def1}) as the definition of the polar rank $\mathrm{prk}(\Gamma)$, valid in any case. 

\begin{note}
{\rm When $\Gamma$ admits singular subspaces of infinite rank, all maximal singular subspaces of $\Gamma$ have infinite rank but not necessarily the same rank (see e.g. \cite{Pas}). Clearly, the polar rank $\mathrm{prk}(\Gamma)$ defined as in (\ref{def1}) is the least upper bound (but possibly not the maximum) of the ranks of the maximal singular subspaces of $\Gamma$.}
\end{note}

A few authors prefer to define the polar rank by means of chains of singular subspaces. Explicitly, let $\mathfrak{C}_{\mathrm{sing}}(\Gamma)$ be the family of all chains of singular subspaces of $\Gamma$. When $\mathrm{prk}(\Gamma)$, defined as in (\ref{def1}), is finite then 
\begin{equation}\label{id2}
\mathrm{sup}(\ell({\cal C})~|~ {\cal C}\in \mathfrak{C}_{\mathrm{sing}}(\Gamma)) ~ = ~ \mathrm{prk}(\Gamma).
\end{equation}
In this case we could take the number 
\[\mathrm{prk}_{\mathrm{C}}(\Gamma) ~ := ~ \mathrm{sup}(\ell({\cal C})~|~ {\cal C}\in \mathfrak{C}_{\mathrm{sing}}(\Gamma))\]
as the polar rank of $\Gamma$, by definition. Johnson \cite{J1, J2} and Cohen \cite{Cohen} indeed define the polar rank in this way, adopting the $\infty$-convention when $\Gamma$ admits singular subspaces of infinite rank. Once again, one might be tempted to extend this definition to the general case, getting rid of the $\infty$-convention, but now this would be an error (even if, in a sense, a definition is never wrong). Indeed, as we know from Theorem \ref{th1}, equation (\ref{id2}) might fail to hold when $\mathrm{prk}(\Gamma)$ is infinite. It certainly fails when $\mathrm{prk}(\Gamma) = \mathrm{max}(\mathrm{rk}(S)~|~ S~\in \mathrm{Sing}(\Gamma))$ and $2^\mathfrak{n} \leq \mathrm{prk}(\Gamma)$ for every $\mathfrak{n} < \mathrm{prk}(\Gamma)$.

Theorem \ref{th2} suggests how to correct the above: instead of considering arbitrary chains of singular subspaces, we must consider only the well ordered ones. Explicitly, let $\mathfrak{W}_{\mathrm{sing}}(\Gamma) \subseteq \mathfrak{C}_{\mathrm{sing}}(\Gamma)$ be the family of all well ordered chains of singular subspaces of $\Gamma$ and put
\[\mathrm{prk}_{\mathrm{WO}}(\Gamma) ~ := ~ \mathrm{sup}(\ell({\cal C})~|~ {\cal C}\in \mathfrak{W}_{\mathrm{sing}}(\Gamma)).\] 
Singular subspaces satisfy (EP), since they are projective spaces. Hence Theorem \ref{th2} can be applied in each of the maximal singular subspaces of $\Gamma$. Thus we obtain the following equality, no matter if $\mathrm{prk}(\Gamma)$ is finite or infinite:  
\begin{equation}\label{id5}
\mathrm{prk}_{\mathrm{WO}}(\Gamma) ~ = ~ \mathrm{prk}(\Gamma).
\end{equation} 
When $\mathrm{prk}(\Gamma)$ if finite then $\mathfrak{C}_{\mathrm{sing}}(\Gamma) = \mathfrak{W}_{\mathrm{sing}}(\Gamma)$ and we get back (\ref{id2}). 

\begin{note}
{\rm In \cite[Introduction]{Pas} it is wrongly claimed that (\ref{id2}) holds in general. However, no mention of the number $\mathrm{prk}_{\mathrm{C}}(\Gamma)$ is made in \cite{Pas} after that claim. So, luckily, that error has no consequences in \cite{Pas}.}
\end{note} 

\begin{note}
{\rm The problem we have discussed in this subsection cannot arise in the setting chosen by Tits \cite{Tits} for polar spaces. Indeed in that setting all polar spaces have finite rank, by definition.}
\end{note} 

\subsection{The generating rank of a polar space}

Throughout this subsection $\Gamma = (P, {\cal L})$ is a non-degenerate polar space of finite polar rank $\mathrm{prk}(\Gamma) \geq 2$, but $\mathrm{rk}_{\mathrm{gen}}(\Gamma)$ is allowed to be infinite. We also assume that $\Gamma$ is {\em thick-lined}, namely all lines of $\Gamma$ have at least three points. 

The main results of this subsection are a remake of Section 2.3 of Cardinali, Giuzzi and Pasini \cite{CGP}. We shall state them in \S\S \ref{sec faith} and \ref{sec corank}, but before to come to them we need to recall a few basics and well known theorems on projective embeddings of polar spaces. 

\subsubsection{Preliminaries on embeddings}

A {\em (projective) embedding} of $\Gamma$ is an injective mapping $e:P\rightarrow \mathrm{PG}(V)$ such that the set $e(\ell) := \{e(x)~|~x\in \ell\}$ is a line of $\mathrm{PG}(V)$ for every line $\ell\in{\cal L}$ and $e(P)$ spans $\mathrm{PG}(V)$. The dimension of the vector space $V$ is taken as the {\em dimension} of the embedding $e$, henceforth denoted by the symbol $\mathrm{dim}(e)$. If $\mathbb{K}$ is the underlying division ring of $V$, we say that $e$ is {\em defined over} $\mathbb{K}$.  

Let $e:P\rightarrow \mathrm{PG}(V)$ be an embedding of $\Gamma$. If $X\in \mathrm{Gen}(\Gamma)$ then $e(X)$ spans $\mathrm{PG}(V)$. Hence
\begin{equation}\label{rk-dim}
\mathrm{rk}_{\mathrm{gen}}(\Gamma) ~\geq~ \mathrm{dim}(e).
\end{equation}
Moreover, if $M, M'$ are disjoint maximal subspaces of $\Gamma$ then $e(M)$ and $e(M')$ are disjoint subspaces of $\mathrm{PG}(V)$ of rank $n = \mathrm{prk}(\Gamma)$. They span a subspace of $\mathrm{PG}(V)$ of rank $2n$. Therefore 
\begin{equation}\label{prk-dim}
\mathrm{dim}(e) ~ \geq ~ 2\cdot\mathrm{prk}(\Gamma).
\end{equation}
By (\ref{rk-dim}) and (\ref{prk-dim}), if $\Gamma$ is embeddable, namely it admits an embedding, then 
\begin{equation}\label{rk-prk}
\mathrm{rk}_{\mathrm{gen}}(\Gamma) ~ \geq ~ 2\cdot\mathrm{prk}(\Gamma).
\end{equation}
The following is well known (Tits \cite[Chapters 8 and 9]{Tits}; see also Buekenhout and Cohen \cite[Chapter 8]{BC}): 

\begin{theo}\label{Tits1}
Suppose that $\mathrm{prk}(\Gamma) \geq 3$. Then $\Gamma$ is embeddable except in the following two exceptional cases, where $\mathrm{prk}(\Gamma) = 3$:
\begin{itemize}
\item[$(1)$] $\Gamma$ is the line-grassmannian of a 3-dimensional projective geometry defined over a non-commutative division ring. 
\item[$(2)$] The singular planes of $\Gamma$ are Moufang but not desarguesian.   
\end{itemize}
\end{theo} 
\begin{note}
{\rm No non-degenerate generalized quadrangle can be generated by three points and, most likely, five points are not enought to generate a polar space as in cases (1) and (2) of Theorem \ref{Tits1}. If so, inequality (\ref{rk-prk}) holds in any case, let $\Gamma$ be embeddable or not.}
\end{note}  

Given two embeddings $e':P\rightarrow \mathrm{PG}(V')$ and $e:P\rightarrow \mathrm{PG}(V)$ of $\Gamma$, a {\em morphism} (an {\em isomorphism}) from $e'$ to $e$ is a morphism (an isomorphism) of projective spaces $f:\mathrm{PG}(V')\rightarrow\mathrm{PG}(V)$ such that $e = f\cdot e'$. If a morphism (an isomorphism) exists from $e'$ to $e$ then we write $e' \rightarrow e$ for short (respectively, $e'\cong e$). Following Tits \cite{Tits}, we say that an emdedding $e$ is {\em dominant} if $e'\rightarrow e$ implies $e'\cong e$. The following theorem is also contained in Tits \cite[\S 8.6]{Tits}: 

\begin{theo}\label{Tits2}
Suppose that $\Gamma$ is embeddable. Then for every embedding $e'$ of $\Gamma$ there exists a dominant embedding $e$ such that $e\rightarrow e'$. Moreover, all dominant embeddings of $\Gamma$ are mutually isomorphic (in free words, $\Gamma$ admits a unique dominant embedding) except in the following two cases, where $\mathrm{prk}(\Gamma) = 2$:  
\begin{itemize}
\item[$(1)$] $\Gamma$ is a grid with lines of size at least $6$. In this case all embeddings of $\Gamma$ are $4$-dimensional and defined over a field. 
\item[$(2)$] $\Gamma$ admits just two non-isomorphic embeddings. They are both $4$-dimensional and defined over the same quaternion algebra. 
\end{itemize}
Moreover, when $\Gamma$ is not as  in cases $(1)$ or $(2)$ and its dominant embedding $e$ is defined over a division ring of characteristic other than $2$,  then $e$ is the unique embedding of $\Gamma$.    
\end{theo}
As a consequence of this theorem, if $\Gamma$ is not a grid then all of its embeddings are defined over the same division ring, say $\mathbb{K}$; in short,  $\Gamma$ is {\em defined over} $\mathbb{K}$.

We complete our survey of embeddings with one more celebrated theorem of Tits \cite[Chapter 8]{Tits}.    

\begin{theo}\label{Tits3}
Let $e:P\rightarrow \mathrm{PG}(V)$ be a dominant embedding of $\Gamma$, with $V$ defined over $\mathbb{K}$. Put $e({\cal L}) := \{e(\ell)~|~\ell \in {\cal L}\}$ and $e(\Gamma) := (e(P), e({\cal L})) \cong \Gamma$. 

If $\mathrm{char}(\mathbb{K}) \neq 2$ then $e(\Gamma)$ is the polar space associated to a non-degenerate reflexive sesquilinear form  $f:V\times V\rightarrow \mathbb{K}$ (see {\rm \cite[\S 8.1]{Tits}}).

If $\mathrm{char}(\mathbb{K}) = 2$ then  $e(\Gamma)$ is the polar space associated to a non-degenerate pseudo-quadratic form defined on $V$ (see {\rm \cite[\S 8.2]{Tits}}). 
\end{theo}

\subsubsection{Faithful embeddings}\label{sec faith}

Let $e:P\rightarrow \mathrm{PG}(V)$ be an emebdding of the polar space $\Gamma = (P, {\cal L})$. Given $X\subseteq P$, keeping the symbol $\langle X \rangle$ to denote the subspace of $\Gamma$ generated by $X$, we denote by $[e(X)]$ the span of $e(X)$ in $\mathrm{PG}(V)$. Clearly, the following inclusion holds for every subspace $S$ of $\Gamma$:
\begin{equation}\label{incl 1}
S ~ \subseteq ~ e^{-1}([e(S)]).
\end{equation} 
In general (\ref{incl 1}) is a strict inclusion, even if $e$ is dominant. For instance, let $X$ be a set of pairwise non-collinear points properly contained in its double perp $X^{\perp\perp}$ (notation as usual for polar spaces). Trivially, $X$ is a subspace of $\Gamma$. However $X  \subset  X^{\perp\perp} =  e^{-1}([e(X)])$. So, if we don't put any restriction on the family of subspaces $S$ to be considered in (\ref{incl 1}), there is no hope to turn (\ref{incl 1}) into an equality.  

Following \cite{CGP}, we say that a subspace of $\Gamma$ is {\em nice} if it contains two mutually disjoint maximal singular subspaces. We say that  the embedding $e$ is {\em faithful} if $S = e^{-1}([e(S)])$ for every nice subspace $S$ of $\Gamma$.  
   
\begin{prop}\label{tight0}
Let $e$ be faithful. Then $\mathrm{dim}(e) = \mathrm{rk}_{\mathrm{gen}}(\Gamma)$.
\end{prop}
{\bf Proof.} This proposition is implicit in the proof of Corollary 3.6 of \cite{CGP}, but we shall give an easier proof here. 

Given two disjoint maximal subspaces $M$ and $M'$ of $\Gamma$, choose bases $A$ and $A'$ of $M$ and $M'$ respectively. Next choose $B\subset P$ such that $e(A)\cup e(A')\cup e(B)$ is a basis of $\mathrm{PG}(V)$. This is possible since $e(A)\cup e(A')$ is independent and $e(P)\supseteq e(A)\cup e(A')$ spans $\mathrm{PG}(V)$. Then $[e(A)\cup e(A')\cup e(B)] = \mathrm{PG}(V)$. However $\langle A\cup A'\cup B\rangle = e^{-1}([e(A)\cup e(A')\cup e(B)]$ because $e$ is faithful.  Hence $\langle A\cup A'\cup B\rangle = P = e^{-1}(\mathrm{PG}(V))$. Therefore $\mathrm{dim}(e) \geq \mathrm{rk}_{\mathrm{gen}}(\Gamma)$. Hence $\mathrm{dim}(e) = \mathrm{rk}_{\mathrm{gen}}(\Gamma)$ by (\ref{rk-dim}).  \hfill $\Box$  

\begin{prop}\label{tight1} 
Every faithful embedding is dominant.
\end{prop}
{\bf Proof.} Let $e':\Gamma\rightarrow\mathrm{PG}(V')$ and $e:\Gamma\rightarrow\mathrm{PG}(V)$ be embeddings of $\Gamma$ and let $f:e'\rightarrow e$ be a morphism but not an isomorphism. Then $f$ is induced by a non injective semilinear mapping $\varphi:V'\rightarrow V$. Let $K$ be the subspace of $\mathrm{PG}(V')$ corresponding to the kernel of $\varphi$. Given any two disjoint maximal singular subspaces $M$ and $M'$ of $\Gamma$, the subspace $[e'(M)\cup e'(M')]$ of $\mathrm{PG}(V')$ spanned $e'(M)\cup e'(M')$ meets $K$ trivially, otherwise at least one of the projective lines $[e(x), e(x')]$ with $x\in M$ and $x'\in M'$ would meet $K$ non trivially, thus forcing $f(e'(x)) = f(e'(x'))$ and consequently $e(x) = e(x')$, which cannot be.  

As $K\cap [e'(M)\cup e'(M')] = \emptyset$, we can choose a hyperplane $H$ of $\mathrm{PG}(V')$ containing $[e'(M)\cup e'(M')]$ but not $K$. Then $S := e'^{-1}(H)$ is a proper subspace of $\Gamma$. Moreover $S$ contains $M\cup M'$, hence it is nice. On the other hand, $f(H) = \mathrm{PG}(V)$. Hence $e^{-1}([e(S)]) = e^{-1}([f(e'(S)]) = e^{-1}(f(H)) = e^{-1}(\mathrm{PG}(V)) = P$. So, $S \subset e^{-1}([e(S)])$. The embedding $e$ is unfaithful.  \hfill $\Box$  

\begin{conj}\label{tight3}
All dominant embeddings are faithful.  
\end{conj} 

The following result by Cardinali, Giuzzi and Pasini \cite[Lemma 2.3]{CGP} partially supports the previous conjecture. 

\begin{prop}\label{tight2}
Suppose that $\Gamma$ is not as in cases (1) or (2) or Theorem {\rm \ref{Tits2}} and let $e$ be its unique dominant embedding. Let $\mathbb{K}$ be the underlying division ring of $\Gamma$ and suppose that either $\mathbb{K}$ is commutative or $\mathrm{char}(\mathbb{K}) \neq 2$. Then the embedding $e$ is faithful. 
\end{prop}

\begin{note}
{\rm Only the commutative case is considered in Lemma 2.3 of \cite{CGP} but the arguments used to prove that Lemma also work when $\mathbb{K}$ is non-commutative provided that, for every nice subspace $S$ of $\Gamma$, the polar space $e(S)$ can be described by a reflexive sesquilinear form defined over (the underlying vector space of) $[e(S)]$. In view of Theorem \ref{Tits3}, this is the case when $\mathrm{char}(\mathbb{K}) \neq 2$.}
\end{note}  

\subsubsection{The polar corank and chains of nice subspaces}\label{sec corank} 

Let $\Gamma$ be embeddable. Then $\mathrm{rk}_{\mathrm{gen}}(\Gamma) \geq 2\cdot\mathrm{prk}(\Gamma)$ by (\ref{rk-prk}). As $\mathrm{prk}(\Gamma)$ is assumed to be finite, there exists a unique cardinal number $\mathfrak{r}$ such that 
\[\mathrm{rk}_{\mathrm{gen}}(\Gamma) ~= ~\mathfrak{r} + 2\cdot\mathrm{prk}(\Gamma).\]
(Clearly, $\mathfrak{r} = \mathrm{rk}_{\mathrm{gen}}(\Gamma)$ when $\mathrm{rk}_{\mathrm{gen}}(\Gamma)$ is infinite.) We call $\mathfrak{r}$ the {\em polar corank} of $\Gamma$, denoting it by the symbol $\mathrm{crk}(\Gamma)$. 

The number $\mathrm{crk}(\Gamma)$ can be characterized independently of $\mathrm{rk}_{\mathrm{gen}}(\Gamma)$ by means of maximal well ordered chains of nice subspaces of $\Gamma$. Explicitly, let $\mathfrak{W}_{\mathrm{nice}}(\Gamma)$ be the family of well ordered chains of nice subspaces of $\Gamma$ and put
\[\mathrm{crk}_{\mathrm{WO}}(\Gamma) ~ := ~ \mathrm{sup}(\ell({\cal C})~|~ {\cal C}\in \mathfrak{W}_{\mathrm{nice}}(\Gamma)).\]
The following theorem generalizes theorems 2.8 and 2.9 of \cite{CGP}:

\begin{theo}\label{ark1}
Suppose that $\Gamma = (P, {\cal L})$ admits a faithful embedding. Then 
\begin{equation}\label{ark0}
\mathrm{crk}(\Gamma) ~ = ~ \mathrm{crk}_{\mathrm{WO}}(\Gamma).
\end{equation}
Moreover, every member of $\mathfrak{W}_{\mathrm{nice}}(\Gamma)$ is contained in a maximal member of $\mathfrak{W}_{\mathrm{nice}}(\Gamma)$ and all maximal members of $\mathfrak{W}_{\mathrm{nice}}(\Gamma)$ have length $\mathrm{crk}(\Gamma)$. 
\end{theo}
{\bf Proof.} Choose a minimal nice subspace $S$. So, $S = \langle M, M'\rangle$ for two disjoint maximal singular subspaces of $\Gamma$ and
$\mathrm{rk}_{\mathrm{gen}}(S) = 2n$, where $n := \mathrm{prk}(\Gamma)$. We firstly prove the following claim:
\begin{itemize}
\item[$(\ast)$] For $S\subseteq X\subseteq P$ and $x, y \in P\setminus S$, suppose that $y\in \langle X\cup\{x\}\rangle$ but $y\not\in\langle X\rangle$. Then $x\in \langle X\cup\{y\}\rangle$.
\end{itemize}
Let $e:P\rightarrow \mathrm{PG}(V)$ be a faithful embedding of $\Gamma$. The subspace $\langle X\rangle$ is nice. Hence $e(\langle X\rangle) = [e(X)]$ because $e$ is faithful. Accordingly, $e(y)\not \in [e(X)]$ because $y\not \in \langle X\rangle$. On the other hand $e(y)\in [e(X)\cup\{e(x)\}]$ because $y\in \langle X\cup\{x\}\rangle$. Consequently $e(x)\in [e(X)\cup\{e(y)\}]$ by (EP) in $\mathrm{PG}(V)$. Therefore $x\in \langle X\cup \{y\}\rangle$ since $e$ is faithful. Claim $(\ast)$ is proved.

By $(\ast)$ we see that, for two points $x, y\in P\setminus S$, we have $y\in \langle S\cup\{x\}\rangle$ if and only if $\langle S\cup\{x\}\rangle = \langle S\cup\{y\}\rangle$. So, we can define a point-line geometry $\Gamma(S)$ by taking the subspaces $\langle S\cup\{x\}\rangle$ with $x\in P\setminus S$ as points and the sets $\langle S\cup\{x,y\}\rangle$ with $\langle S\cup\{x\}\rangle\neq \langle S\cup\{y\}\rangle$ as lines. The points of a line $L = \langle S\cup\{x,y\}\rangle$ are the subspaces $\langle S\cup\{z\}\rangle$ for $z\in L\setminus S$. By claim $(\ast)$, the Exchange Property (EP) holds in $\Gamma(S)$. Therefore 
\begin{equation}\label{ark2}
\mathrm{rk}_{\mathrm{gen}}(\Gamma(S)) ~ = ~ \mathrm{rk}_{\mathrm{WO}}(\Gamma(S))
\end{equation}
by Theorem \ref{th2}. We shall now prove the following:
\begin{equation}\label{claim2}
\mathrm{rk}_{\mathrm{gen}}(\Gamma(S)) ~ = ~ \mathrm{crk}(\Gamma).
\end{equation}
Let $B$ and $B'$ be bases of $M$ and $M'$ respectively. Then, for every subset $X\subset S\setminus P$, the set $B\cup B'\cup X$ generates $\Gamma$ if and only if the set $S[X] := \{\langle S\cup\{x\}\rangle\}_{x\in X}$ generates $\Gamma(S)$. Moreover, if $\langle S\cup\{x\}\rangle = \langle S\cup \{y\}\rangle$ for two distinct points $x, y \in X$ then $\langle (B\cup B'\cup X)\setminus \{y\}\rangle = \langle B\cup B'\cup X\rangle$. It follows that
\begin{equation}\label{ark3}
\mathrm{rk}_{\mathrm{gen}}(\Gamma) ~\leq~  2n + \mathrm{rk}_{\mathrm{gen}}(\Gamma(S)).
\end{equation}
On the other hand, we can always choose $X$ in such a way that $e(B)\cup e(B')\cup e(X)$ is a basis of $\mathrm{PG}(V)$. Clearly, with $X$ chosen in this way, $S[X]$ generates $\Gamma(S)$. Therefore
\begin{equation}\label{ark4}
\mathrm{dim}(e) ~\geq~  2n + \mathrm{rk}_{\mathrm{gen}}(\Gamma(S)).
\end{equation}
By (\ref{ark3}), (\ref{ark4}) and Proposition \ref{tight0} we obtain that $\mathrm{rk}_{\mathrm{gen}}(\Gamma) =  2n + \mathrm{rk}_{\mathrm{gen}}(\Gamma(S))$. Hence $\mathrm{rk}_{\mathrm{gen}}(\Gamma(S)) = \mathrm{crk}(\Gamma)$, as claimed.

By (\ref{claim2}) and (\ref{ark2}) we obtain that $\mathrm{crk}(\Gamma) = \mathrm{rk}_{\mathrm{WO}}(\Gamma(S))$ for every minimal nice subspace $S$. Equality (\ref{ark0}) follows. The remanining claims of the theorem follow from Corollary \ref{cor 3.6} and Theorem \ref{theo 3.7} applied to $\Gamma(S)$, with $S$ any minimal nice subspace.  \hfill $\Box$  

\begin{cor}
Let $e:P\rightarrow \mathrm{PG}(V)$ be a faithful embedding of $\Gamma$ and let $f$ be a non-degenerate reflexive sesquilinear form of $V$ defining $e(\Gamma)$ or the sequilinearization of a non-degenerate pseudoquadratic form defining $e(\Gamma)$. 

Given two mutually disjoint maximal singular subspaces $M$ and $M'$ of $\Gamma$, let $[e(M)\cup e(M')]^\perp$ be the orthogonal complement of $[e(M)\cup e(M')]$ in $\mathrm{PG}(V)$, orthogonality being defined with respect to $f$. Then 
\[\mathrm{crk}_{\mathrm{WO}}(\Gamma) = \mathrm{rk}([e(M)\cup e(M')]^\perp).\] 
\end{cor} 

\begin{note}
{\rm  The polar corank $\mathrm{crk}(\Gamma)$ is called {\em anisotropic defect} in \cite{CGP}, in view of the fact that $[e(M)\cup e(M')]^\perp\cap e(P) = \emptyset$.  However the word `defect' is usually given a different meaning in the literature. This considered, we have replaced the name `anisotropic defect' with `polar corank', which hopefully transmits no wrong suggestions.}
\end{note}

\vspace{10 mm}

\noindent
Author's address:\\

\noindent
Antonio Pasini\\
Department of Information Engineering and Maths.\\
University of Siena\\
Via Roma 56, ~ 53100 Siena, Italy\\
antonio.pasini@unisi.it;  antoniopasini13@gmail.com

\end{document}